\documentclass[11pt]{amsart}

\usepackage{amssymb}
\usepackage[hidelinks]{hyperref}

\title[Projective surfaces and profinite completion]{Projective surfaces
  with infinite fundamental group and trivial profinite completion}

\author{Wenchuan Hu}
\author{Peihan Wu}
\address{School of Mathematics, Sichuan University, Chengdu 610064, P. R. China}
\email[Wenchuan Hu]{huwenchuan@gmail.com, wenchuan@scu.edu.cn}
\email[Peihan Wu]{3119567113@qq.com}

\hypersetup{
  pdftitle={Projective surfaces with infinite fundamental group and trivial profinite completion},
  pdfauthor={Wenchuan Hu and Peihan Wu}
}

\theoremstyle{plain}
\newtheorem{theorem}{Theorem}
\newtheorem{proposition}[theorem]{Proposition}
\newtheorem{lemma}[theorem]{Lemma}
\newtheorem{corollary}[theorem]{Corollary}
\DeclareMathOperator{\SO}{SO}
\DeclareMathOperator{\CAT}{CAT}
\newcommand{\RR}{\mathbb R}
\newcommand{\CC}{\mathbb C}
\newcommand{\ZZ}{\mathbb Z}
\newcommand{\QQ}{\mathbb Q}
\newcommand{\PP}{\mathbb P}
\newcommand{\nc}[1]{\langle\!\langle #1\rangle\!\rangle}
\newcommand{\T}{\mathcal T}
\newcommand{\et}{\mathrm{\acute et}}
\begin{document}

\begin{abstract}
For every nontrivial finite abelian group \(A\) of odd order, we construct
a smooth complex projective surface with infinite topological fundamental
group and trivial \'{e}tale fundamental group. The isomorphism types of its
finite subgroups are precisely those of the subgroups of \(A\).
The surface is Kobayashi hyperbolic and has Stein universal cover.
\end{abstract}

\maketitle
%\tableofcontents

\section{Introduction}
Riemann existence identifies the \'etale fundamental group of a smooth complex projective variety with the profinite completion of its topological fundamental group \cite[Exp.~XII, Th.~5.1]{SGA1}. Koll\'ar \cite[p.~70]{Kollar} and Poonen \cite[Rem.~3.5.47]{Poonen} ask whether the former can be trivial while the latter is nontrivial. Hughes, Llosa Isenrich, Py, Stover and Vidussi ask whether every infinite K\"ahler group has a nontrivial finite quotient \cite[\S~1; arXiv v2, p.~2]{Hughes}. We construct surfaces that answer these questions.

Toledo gave projective examples with non-residually finite fundamental group \cite{Toledo}. Eyssidieux subsequently used root stacks of totally geodesic divisors to produce odd-order inertia that survives in the orbifold fundamental group and dies in every finite image \cite[Prop.~1.10]{Eys}. Those groups map onto the original arithmetic lattice. We start instead with a quotient by a reflection. We replace its weight two by an odd weight along every mirror. Killing all mirror meridians in the resulting orbifold group leaves a finite group, and the arithmetic boundary bound forces every finite image to factor through it.

\textbf{Disclosure of AI assistance.}
The authors used Open AI Codex  for assistance with manuscript preparation
and revision. They take responsibility for the text and the
mathematical arguments.

\section{Main results}
We take the fundamental group of an analytic orbifold to mean the fundamental group of its classifying space. We denote normal closure by \(\nc{\cdot}\) and the cyclic group of order \(d\) by \(C_d\).

\begin{theorem}\label{thm:main}
Let \(A\) be a nontrivial finite abelian group of odd order. There is a smooth connected complex projective surface \(X\) with ample canonical bundle, infinite fundamental group \(K=\pi_1(X)\) and \(\pi_1^{\et}(X)=1\). The finite subgroups of \(K\), up to isomorphism, are precisely the subgroups of \(A\). The surface \(X\) is Kobayashi hyperbolic and its universal cover is Stein. For every irreducible curve \(C\subset X\) with normalization \(C^\nu\), the image of \(\pi_1(C^\nu)\to K\) is infinite.
\end{theorem}

\begin{corollary}\label{cor:applications}
For the groups \(K\) of Theorem~\ref{thm:main}, every finite-dimensional linear representation over any field is trivial. For every prime \(p\), the degree-two comparison with trivial \(\mathbb F_p\) coefficients has zero source \(H^2_{\mathrm{cont}}(\widehat K,\mathbb F_p)\) and nonzero target \(H^2(K,\mathbb F_p)\). In particular \(K\) fails Serre goodness.
\end{corollary}

\section{The arithmetic quotient and its mirrors}
Let \(k=\QQ(\sqrt2)\), let \(R=\ZZ[\sqrt2]\), and put
\[
 q(x)=\sqrt2(x_1^2+x_2^2)-x_3^2-x_4^2-x_5^2-x_6^2 .
\]
Write \(G_0=\SO(q,\RR)^\circ\), and let \(S\) be a component of the space of oriented positive two-planes. It is the type-IV Hermitian symmetric domain with its complete invariant K\"ahler metric. Set
\[
 \Gamma=G_0\cap\ker\bigl(\mathrm O(q,R)\longrightarrow\mathrm{GL}_6(R/3R)\bigr).
\]
At the conjugate real place \(q\) is negative definite. The arithmetic lattice and compactness criteria therefore make \(\Gamma\) a uniform lattice in \(G_0\) \cite{BHC,Toledo}. The level-three kernel is torsion-free by Minkowski's congruence lemma, applied after restricting scalars to \(\ZZ\). In particular \(M=\Gamma\backslash S\) is smooth. It is projective by positivity of the canonical bundle \cite[Lem.~2.2]{Toledo}.

The reflection \(r\) negating \(e_6\) normalizes \(\Gamma\) and acts holomorphically on \(S\). Write \(G=\Gamma\rtimes\langle r\rangle\) and \(Y=M/\langle r\rangle\). The latter is projective. If \(s=\gamma r\) fixes a point of \(S\), then \(s^2=1\) by freeness of \(\Gamma\). Its negative eigenspace has odd dimension \(b\), since \(\det s=-1\), and
\[
 6-2b=\operatorname{tr}(s)\equiv \operatorname{tr}(r)=4\pmod{3R}.
\]
It follows that \(b=1\). The invariant oriented positive plane is fixed pointwise, so the negative eigenline is orthogonal to it. Thus every such \(s\) is a reflection in a totally geodesic complex hypersurface of type \(\SO^+(2,3)\). The fixed locus of the descended involution on \(M\) consists of finitely many disjoint smooth compact divisors \(D_i\). Their images \(\bar D_i\) are disjoint smooth Cartier divisors in the smooth fourfold \(Y\), with coordinates \(t=u^2\) normal to \(\bar D_i\).

\begin{lemma}\label{lem:boundary}
For every \(D_i\), the square of a normal circle meridian of \(D_i\) vanishes in every finite image of the fundamental group of its normal circle bundle.
\end{lemma}

\begin{proof}
Choose a mirror \(D_S\) over \(D_i\) and let \(\Gamma_i\) be its stabilizer in \(\Gamma\). It preserves the negative eigenline of its reflection and fixes a vector \(v\in k^6\) on that line. Indeed the only possible signs are \(\pm1\), and \(\gamma+I\) is invertible over \(k\) for \(\gamma\in\Gamma\) since its determinant is nonzero modulo three.

The restriction of \(\Gamma_i\) to the mirror is a uniform lattice in \(H=\SO^+(2,3)\). Indeed distinct translates are disjoint, since two reflections fixing the same point of \(S\) have product in \(\Gamma\), and \(\Gamma_i\backslash D_S=D_i\) is compact. The normal circle bundle is \(\Gamma_i\backslash H/\SO(3)\), by the homogeneous normal-plane description of \cite[Lem.~3.3]{Toledo}. Its fundamental group is the preimage lattice in the universal covering group of \(H\), modulo the \(\ZZ/2\) from \(\pi_1(\SO(3))\). The covering kernel over \(\mathrm{Spin}(2,3)^\circ\) is the graph of reduction modulo two in \(\pi_1(H)=\ZZ\oplus\ZZ/2\). It projects isomorphically onto the meridian group \cite[Lem.~3.4]{Toledo}.

Set \(q_v=q|_{v^\perp}\) and \(G_v=\mathrm{Spin}(q_v)\). The restriction of \(\Gamma_i\) is arithmetic in \(\SO(q_v)\) because \(\Gamma_i\) is the intersection of \(\Gamma\) with the \(k\)-subgroup fixing \(v\). By restriction of scalars and the arithmetic isogeny theorem \cite[Cor.~6.11]{BHC}, the image of a spin arithmetic subgroup is commensurable with this restriction. At the two real places \(q_v\) has Witt indices \(2\) and \(0\). The latter implies anisotropy over \(k\). Since \(\dim q_v=5\) and the total archimedean rank is two, Kneser's theorem gives perfection of \(G_v(k)\) and centrality of its archimedean congruence kernel \cite[p.~406]{Rag}. Writing \(k_\infty=k\otimes_{\QQ}\RR\), we have \(G_v(k_\infty)\simeq\mathrm{Sp}_4(\RR)\times\mathrm{Spin}(5)\), and the compact factor is simply connected. Deligne's comparison puts the kernel of the restriction pairing
\[
 \pi_1(G_v(k_\infty))\longrightarrow M(\varnothing,G_v)^\vee
\]
in every finite-index subgroup of the full inverse image of a spin arithmetic subgroup \cite[Th.~1.4]{Rag}. Here \(M(\varnothing,G_v)\) is the absolute metaplectic kernel. Its order is at most two by the calculation for type \(B_2=C_2\) in \cite[Th.~8.1, p.~167]{PR}, since \(k\) has only the roots of unity \(\{\pm1\}\). Given a finite-index subgroup of the circle-bundle group, take its inverse image in the universal covering group of \(H\) and intersect with the lifted spin arithmetic subgroup. The intersection has finite index in the spin lift by the commensurability above. Thus it contains twice the cyclic covering kernel, which projects to twice the meridian by the preceding identification.
\end{proof}

For the odd-weight construction, Toledo's original \(8\ZZ\) bound \cite[Lem.~3.4]{Toledo} also suffices. That bound uses the main theorem of \cite{Rag}, whose local splitting property is discussed in \cite{RagCorr}. Lemma~\ref{lem:boundary} uses the metaplectic calculation and Deligne's comparison directly.

Put \(U=Y\setminus\bigcup_i\bar D_i\), \(\Pi=\pi_1(U)\), and let \(\mu_i\) be a meridian of \(\bar D_i\). On \([M/\langle r\rangle]\) each pullback of \(\bar D_i\) is twice its reduced inverse image as a Cartier divisor. The corresponding line bundle with section defines a map to the square-root stack, which is an isomorphism in the charts \(t=u^2\) and off the mirrors. Write \(\T_e\) for the root stack of weight \(e\) along all \(\bar D_i\), with \(\T_1=Y\) \cite[Ex.~2.4.1]{Cadman}.

For the meridian presentation, take disjoint tubular neighborhoods of the \(\bar D_i\). Over \(\bar D_i\) the classifying space of the root neighborhood retracts onto the classifying space of the root divisor, with homotopy fiber \(BC_e\) over \(\bar D_i\). The punctured neighborhood retracts onto the normal circle bundle, and the map on fiber fundamental groups is \(\ZZ\to C_e\), reduction modulo \(e\), by the local model \(t=w^e\). The two homotopy exact sequences therefore show that the punctured neighborhood surjects on the fundamental group of the root neighborhood, with kernel generated by \(\mu_i^e\). Van Kampen for the classifying spaces gives \(\pi_1(\T_e)=\Pi/\nc{\mu_i^e:\text{ all }i}\), as in \cite[\S~1.1]{Eys}. In particular, for odd \(d\ge3\) we have
\begin{equation}\label{eq:meridian}
 L_d:=\pi_1(\T_d)=\Pi/\nc{\mu_i^d:\text{ all }i},\qquad
 G=\Pi/\nc{\mu_i^2:\text{ all }i}.
\end{equation}
Here the second identification also uses contractibility of \(S\). Map the normal circle bundle of \(D_i\) into \(M\setminus\bigcup_j D_j\) along the boundary of a small tube, and then to \(U\). Its fiber meridian maps to \(\mu_i^2\) up to conjugacy. Restricting any finite image of \(L_d\) along these maps and \(\Pi\to L_d\), Lemma~\ref{lem:boundary} kills \(\mu_i^4\). Using Toledo's \(8\ZZ\) bound instead would kill \(\mu_i^{16}\). Either relation, together with \(\mu_i^d=1\) for odd \(d\), kills \(\mu_i\) in the finite image. Thus every finite image factors through
\begin{equation}\label{eq:coarse}
 F=\Pi/\nc{\mu_i:\text{ all }i}=\pi_1(Y).
\end{equation}

The group \(F\) is finite. The meridian of the image of the \(r\)-mirror maps to \(r\) up to conjugacy in \(G\), so \(F\) is a quotient of \(G/\nc{r}\). If \(\Gamma\) centralized \(r\), the distance to the mirror would descend to the compact manifold \(M\). This distance is unbounded along a normal geodesic ray. Hence \(\Gamma\cap\nc{r}_G\) contains a nontrivial commutator and is infinite by torsion-freeness. The group \(G_0\) has simple Lie algebra, finite center and real rank two. Passing to its adjoint group and applying the Margulis normal subgroup theorem shows that \(\Gamma\cap\nc{r}_G\) has finite index \cite[Th.~VIII.2.6]{Margulis}. The image of \(\Gamma\) generates \(G/\nc{r}\), proving the claim.

\section{The modified orbifold and its development}
The metric on \(Y\) is the quotient length metric from \(S\). On the complement of the mirror divisors pull it back to the orbifold charts of \(\T_d\) and complete locally at the divisors.

\begin{proposition}\label{prop:metric}
The completed chart metrics are locally \(\CAT(0)\) and have the given coarse metric as their metric quotient.
\end{proposition}
\begin{proof}
Work in \(S\) with one mirror \(D_S\). The normal exponential coordinates for this closed convex totally geodesic submanifold identify \(S\) with its normal bundle, a trivial oriented plane bundle over the contractible domain \(D_S\). Complete the connected degree-\(d\) cyclic cover of \(S\setminus D_S\) in its pulled-back length metric to obtain \(Z\). There is one point above each point of the mirror. This follows from the local comparison in normal coordinates with a product of a tangential ball and a cone of angle \(2d\pi\). Projection preserves distances to points over the mirror, since the base geodesic to such a point lifts with its length preserved. Thus \(Z\to S\) is a simple branched cover in the sense of \cite[\S~2, pp.~443--444]{All99}. Since \(S\) is complete \(\CAT(0)\) and \(D_S\) is closed and convex, \(Z\) is \(\CAT(0)\) by \cite[Th.~2.1, p.~444]{All99}; it is complete by construction.

On \(Z\) the angular shift by \(d\pi\) is an isometric involution \(\tau\) lifting the reflection, since \(d\) is odd. Write \(W=Z/\langle\tau\rangle\). Its metric is the minimum of the two distances between lifts. In particular it is complete and geodesic. The projection is \(1\)-Lipschitz and preserves distances whenever one endpoint is fixed. The mirror is convex in \(W\), since equality in the triangle inequality between fixed points lifts to equality in \(Z\).

We check the hypotheses of \cite[Th.~3.1]{All13}. For a geodesic triangle with a mirror vertex, lift the opposite side first and then the sides to the unique lift of the mirror vertex. A side lifts by path lifting off the mirror and by the distance identity at its endpoints on the mirror. Its intersection with the mirror is an interval, by convexity. Each lift is geodesic because it preserves length and the projection is \(1\)-Lipschitz. Comparison follows in \(Z\). By \cite[Lem.~5.2]{All13}, the tangent of \(Z\) at the mirror is \(\RR^6\) times the Euclidean cone on a circle of length \(2d\pi\). Geodesic germs based at a fixed point project to geodesic germs by the distance identity, and every germ below lifts as above. The distance minimum thus passes to the tangent by taking limits. The tangent of \(W\) is the same product with circle length \(d\pi\). It is complete \(\CAT(0)\) and has unique local geodesics between given endpoints. Finally, at distance \(h>0\) from the mirror, the two lifts \(z,\tau z\) are at distance \(2h\). Their midpoint lies on the mirror, and a path through a nearest mirror point gives the upper bound. The closed \(h/4\)-ball about \(z\) maps isometrically onto the corresponding ball in \(W\), since within that ball distances between opposite lifts are at least \(3h/2\). The ball is complete \(\CAT(0)\). This proves that \(W\) is \(\CAT(0)\).

In holomorphic coordinates the cover is \(v^d=u\). The involution sends \(v\) to \(-v\), so its quotient coordinate \(w=v^2\) satisfies \(w^d=u^2=t\). Thus \(W\) gives the local root chart. The deck rotations of \(Z\to S\) and \(\tau\) generate the rotations by \(\pi\), so the order-\(d\) quotient of \(W\) is the original reflection quotient. These are also metric quotients, by lifting base paths and taking the minimum over lifts. This identifies the chart metric with the completed pullback of the coarse metric.
\end{proof}

Take small convex balls in \(W\) about fixed mirror points as charts. Each is the inverse image of the corresponding small coarse ball because distances to its center are unchanged by the finite quotient. Away from the mirrors use regular metric charts. Chart metrics agree on overlaps because they agree off the mirrors and are locally completed length metrics. On a sufficiently small overlap, a path leaving a larger overlap has a positive projected length cost, so it does not change distances in the germ. The chart topology is the analytic topology by the local cone comparison, and small closed chart balls are compact. The atlas groupoid is Hausdorff and its source-target map is proper with finite fibers by the separated diagonal. A continuous arrow path with convergent source and target has precompact tails by properness. The intersection of the nested connected closures of the tails is a connected subset of a finite fiber, hence a singleton. Thus the arrow path converges, as required in \cite[III.G.2.10]{BH}. The induced orbit length metric is the compact coarse metric by Proposition~\ref{prop:metric} and subdivision of coarse paths into local charts. The metric Cartan--Hadamard theorem for groupoids \cite[III.G.2.15]{BH} gives a simply connected complete \(\CAT(0)\) development \(B_d\). Its complex manifold structure comes from the root charts. The local charts are locally compact, so \(B_d\) is proper as a complete locally compact length space.

We have \(\T_d^{\mathrm{an}}=[B_d/L_d]\), with point stabilizers equal to the chart inertia groups. Thus every meridian of a mirror has order exactly \(d\). Every finite subgroup of \(L_d\) fixes a point of \(B_d\) by \cite[II.2.8]{BH}, so it is contained in a conjugate of a mirror inertia group unless it is trivial. Put \(K_d=\ker(L_d\to F)\). This group contains all the inertia groups. If a homomorphism from \(K_d\) to a finite group has kernel \(J\), the core of \(J\) in \(L_d\) has finite index in \(L_d\). The quotient by the core factors through \(F\), so the core contains \(K_d\) and \(J=K_d\). Consequently \(\widehat K_d=1\). The \(r\)-mirror gives a copy of \(C_d\) in \(K_d\), so \(K_d\) is nontrivial; since it has no nontrivial finite quotient, it is infinite. Its finite subgroups have precisely the isomorphism types of subgroups of \(C_d\).

\begin{proposition}\label{prop:stein}
The complex manifold \(B_d\) is Stein.
\end{proposition}
\begin{proof}
Fix \(b\in B_d\) and put \(f(x)=\operatorname{dist}(b,x)^2\). Along a unit-speed local geodesic \(\gamma\) the \(\CAT(0)\) inequality gives convexity of \(f(\gamma(t))-t^2\). Choose local potentials \(\phi\) for the original invariant K\"ahler metric, invariant under the reflection in each mirror chart. Finitely many relatively compact subcharts cover the orbifold, so a single \(c>0\) satisfies \(c\operatorname{Hess}_g\phi\le g\) there. On the regular part let \(\Phi\) be the pullback of \(\phi\) to the root chart. Then \(f-c\Phi\) is locally strongly geodesically convex. By smooth approximation \cite[Prop.~2.2]{GreeneWu}, it is plurisubharmonic, since smooth geodesically convex functions on a K\"ahler manifold have nonnegative Levi form.

Near a mirror, \(\Phi(a,w)=\phi(a,u)\) with \(u^2=w^d\) is single-valued by reflection invariance. Evenness gives \(\phi(a,u)-\phi(a,0)=O(|u|^2)\) and first normal derivatives \(O(|u|)\). Local branches of \(u\) have derivatives of order \(j\le2\) bounded by \(O(|w|^{d/2-j})\). Thus the substituted difference and its first and second normal derivatives are respectively \(O(|w|^d)\), \(O(|w|^{d-1})\), and \(O(|w|^{d-2})\), with the same estimates for tangential differentiations up to total order two. This gives the \(C^2\) extension of \(\Phi\) and zero mixed and normal Hessian terms on the mirror. The continuous \(f-c\Phi\) is plurisubharmonic across the mirror by removability. The forms \(dd^c\Phi\) agree off the mirrors and extend continuously to a global semipositive form \(\omega_d\) on \(\T_d\), whose restriction to a mirror is its original K\"ahler form and whose kernel there is the normal direction.

Let \(\psi\) be the squared Hermitian norm of the defining section of the reduced mirror divisor on \(\T_d\). This is a bounded smooth function whose Levi form at the mirror is positive on the normal direction and zero on tangential vectors and mixed pairs. The form \(c\omega_d+dd^c\psi\) is positive on a neighborhood of the mirrors. By convex combination with \(c\omega_d\), it remains positive there with \(\epsilon dd^c\psi\) for \(0<\epsilon\le1\). On the compact complement of that neighborhood, \(\omega_d\) is positive, so a single sufficiently small \(\epsilon>0\) works everywhere. Lifting \(\psi\) to \(B_d\), we obtain the continuous strictly plurisubharmonic function \(f+\epsilon\psi\) since it is locally the sum of \(f-c\Phi\) and a strictly plurisubharmonic \(C^2\) function. It is an exhaustion because \(f\) is proper and \(\psi\) is bounded. The Levi problem \cite[Cor.~5.3.2]{FN} gives the assertion.
\end{proof}

\section{Projective realization and proof of the main results}
\begin{lemma}\label{lem:realization}
Let \(\T\) be a smooth connected effective proper complex Deligne--Mumford stack with projective coarse space \(Y_\T\). There is a smooth connected projective surface \(V\) with ample canonical bundle and a map \(V\to\T\) inducing an isomorphism on fundamental groups. If \(\dim Y_\T=4\), the map can be chosen so that \(V\to Y_\T\) is finite.
\end{lemma}
\begin{proof}
We use the projective-bundle construction of \cite[\S~7.1, Lem.~7.1, pp.~35--36]{OAI} and give the proof needed here. Write \(n=\dim Y_\T\), choose \(a\ge3\), and put
\[
 E=(\mathcal O_\T\oplus T_\T)^{\oplus a}\oplus\mathcal O_\T^{\oplus3},
 \qquad \mathcal P=\PP_\T(E).
\]
Here projectivization parametrizes lines.

At a geometric point \(x\) of \(\T\), let \(g\ne1\) be a stabilizer element. Finite holomorphic actions are linearizable at a fixed point, so effectiveness gives a faithful action on \(T_x\T\). The trivial summand then ensures that \(g\) has at least two eigenvalues on \(\mathcal O_x\oplus T_x\T\). Each eigenspace has codimension at least \(a\) on \(E_x\), since it omits at least one dimension from each of the \(a\) copies. The fixed points of \(g\) on the projective fiber are the projectivizations of those eigenspaces. In finite quotient charts this gives codimension at least \(a\) for the locus of nontrivial stabilizers in \(\mathcal P\), and the same estimate on its image in the coarse space \(P\). Let \(P^\circ\) be the open set coming from trivial stabilizers. It is smooth and identifies with the corresponding open substack of \(\mathcal P\).

The coarse morphism \(\mathcal P\to P\) is a good moduli space morphism: in local finite quotient charts in characteristic zero, taking invariants is exact. Stabilizer orders are bounded, so a sufficiently divisible power of \(\mathcal O_{\mathcal P}(1)\) has trivial stabilizer action on its fibers and descends to \(P\) \cite[Th.~10.3]{Alper}. On each reduced geometric fiber of \(P\to Y_\T\), the pullback of this line bundle under the finite quotient from the corresponding projective space is ample. Thus the descended line bundle is relatively ample. Tensoring with a sufficiently ample pullback from \(Y_\T\) gives an ample line bundle on \(P\), so \(P\) is projective.

For the fundamental group, let \(F_\T\) be the frame bundle of \(T_\T\) and \(Q=\mathcal P\times_\T F_\T\). Faithfulness of the tangent actions makes \(F_\T\) a manifold analytically, and \(Q\to F_\T\) is a projective-space bundle. Its fibers are simply connected, so \(\pi_1(Q)\to\pi_1(F_\T)\) is an isomorphism. The inverse image \(Q^\circ\) of \(P^\circ\) in \(Q\) omits a complex analytic subset of codimension at least \(a\). General position of loops and homotopies gives \(\pi_1(Q^\circ)\simeq\pi_1(Q)\). The homotopy exact sequences of the principal \(\mathrm{GL}_n(\CC)\)-bundles over \(P^\circ,\mathcal P,\T\), taken on classifying spaces, identify their fundamental groups with the corresponding quotients of \(\pi_1(Q^\circ),\pi_1(Q),\pi_1(F_\T)\) by the images of \(\pi_1(\mathrm{GL}_n(\CC))\). Naturality yields
\[
 \pi_1(P^\circ)\simeq\pi_1(\mathcal P)\simeq\pi_1(\T).
\]

The fibers of \(P\to Y_\T\) are finite quotients of projective spaces of dimension \(r_0=\operatorname{rk}(E)-1\), so \(\dim P=n+r_0\). Choose a very ample class \(H\) and cut by \(c=n+r_0-2\) general hyperplanes for the embedding by \(|mH|\). Because \(P\setminus P^\circ\) has codimension at least \(a>2\), the resulting surface \(V\) lies in \(P^\circ\), where it is smooth by Bertini. At each cut on the ordinary locus, the quasi-projective Lefschetz theorem gives an isomorphism on \(\pi_1\) and preserves connectedness \cite[Th.~1.1.3(ii), p.~126]{HammLe}. The last cut starts in complex dimension three. Thus \(V\to\T\) induces an isomorphism on \(\pi_1\). Since \(P\) has quotient singularities, \(K_P\) is \(\QQ\)-Cartier. For large \(m\), \(K_P+cmH\) is ample, and adjunction on \(P^\circ\) gives \(K_V=(K_P+cmH)|_V\).

Suppose now that \(n=4\), so \(c=r_0+2\), and take \(m\ge2\). For any three distinct points of \(P\), products of linear forms in the original embedding yield sections of \(mH\) that are nonzero at one point and vanish at the other two. Thus evaluations at the three points are independent. The ordered triples of distinct points in a common fiber of \(P\to Y_\T\) have dimension at most \(4+3r_0\). A fixed triple imposes three independent linear conditions on each of the \(c\) sections. The incidence of such triples and cutting tuples has dimension at most the dimension of the tuple space plus \(4+3r_0-3c\), which is the dimension of the tuple space minus two. The closure of its projection is a proper closed subset. We can choose the successive general cuts above with their tuple outside it. Every fiber of \(V\to Y_\T\) then has at most two distinct points. The map is proper and quasi-finite, hence finite.
\end{proof}

\subsection*{The cyclic case}
Apply Lemma~\ref{lem:realization} to \(\T_d\), and let \(X_d\to V\) be the connected finite \'etale cover corresponding to \(K_d\subset L_d\), by Riemann existence \cite[Exp.~XII, Th.~5.1]{SGA1}. Its canonical bundle is ample and its fundamental group is \(K_d\). Thus its \'etale fundamental group is trivial and it has the required finite subgroup types.

The analytic fiber product \(X_d\times_{\T_d^{\mathrm{an}}}B_d\), with \(X_d\) viewed analytically, is the pullback of the orbifold universal cover. Each connected component is the universal cover of \(X_d\), since \(\pi_1(X_d)\to L_d\) is injective. The map of this fiber product to \(X_d\times_{Y^{\mathrm{an}}}B_d\) is finite by base change of the relative diagonal of \(\T_d\) over \(Y\). That diagonal is proper and unramified, hence finite. The projection from \(X_d\times_{Y^{\mathrm{an}}}B_d\) to \(B_d\) is finite by Lemma~\ref{lem:realization}. Hence the universal cover of \(X_d\) is Stein by Proposition~\ref{prop:stein} and the finite mapping theorem for Stein spaces \cite[Ch.~V, \S~1.1, Th.~1(d), pp.~125--126]{GR}.

If the image of the fundamental group of a normalized curve in \(X_d\) were finite, the connected finite cover corresponding to the kernel would give a nonconstant holomorphic map from a compact Riemann surface into the Stein universal cover. This is impossible by holomorphic separability. Thus every curve has infinite fundamental-group image.

If \(X_d\) were not Kobayashi hyperbolic, Brody's theorem \cite{Brody} would give a nonconstant entire curve with bounded derivative for a fixed Hermitian metric on \(X_d\). Its coarse image is nonconstant by finiteness over \(Y\). The pullback of the continuous form \(\omega_d\) to \(X_d\) is bounded by a multiple of this Hermitian metric. Hence the density \(\rho\) of the metric \(\rho^2|dz|^2\) on the entire curve is bounded. The invariant metric on \(S\) has holomorphic sectional curvature at most \(-\kappa<0\) by \cite[Prop.~1.3]{DeFab}. Off the discrete set of critical points and mirror contacts, the Gauss equation gives \(\Delta\log\rho\ge\kappa\rho^2\). If the coarse image lies in a mirror, use instead its totally geodesic metric and omit only the critical points. The locally upper-bounded \(\log\rho\) extends subharmonically across the exceptional points. It is bounded above on \(\CC\), hence constant, contradicting the strict inequality where \(\rho>0\). This proves hyperbolicity.

\subsection*{Products}
Write \(A=\prod_i C_{d_i}\) with \(d_i\ge3\) odd. The product of the \(X_{d_i}\) has fundamental group \(\prod_i K_{d_i}\). Every finite subgroup embeds into the product of its projections, and each projection embeds into \(C_{d_i}\). Conversely the group contains a copy of \(A\), so all subgroup types of \(A\) occur. The product has trivial profinite completion since each factor does. If its dimension is greater than two, take successive general ample sections to obtain a smooth surface \(X\) with the same fundamental group by Lefschetz. Its canonical bundle is ample by adjunction. It is hyperbolic as a submanifold of a hyperbolic product. Its universal cover is the connected preimage of \(X\) in the product of the Stein universal covers, since \(X\) induces an isomorphism on fundamental groups. This cover is Stein because finite products of Stein spaces and their closed complex subspaces are Stein \cite[Ch.~V, \S~1.1, Th.~1(e) and (b), pp.~125--126]{GR}. The curve argument then applies to \(X\), completing the proof of Theorem~\ref{thm:main}.

\section{Further consequences}
The image of any finite-dimensional representation of \(K\) over a field is finitely generated and linear, hence residually finite by Mal'cev's theorem \cite[Ths.~VII--VIII]{Malcev}. Since \(\widehat K=1\), this image is trivial. In particular the surfaces answer the question of Eyssidieux and Funar about infinite projective surface groups admitting no complex finite-dimensional representation of infinite image \cite[\S1, p.~2]{EysFun}.

For the cohomological assertion it is enough to consider \(K_d\), since each factor projection \(\prod_i K_{d_i}\to K_{d_j}\) is split. The classifying space of \(\T_d^{\mathrm{an}}\) is a \(K(L_d,1)\) by contractibility of \(B_d\). The coarse map \(X_d\to Y\) therefore factors up to homotopy through the classifying map \(X_d\to BK_d\). For an ample line bundle on \(Y\), its rational first Chern class pulls back to a class of positive self-intersection on \(X_d\), since \(X_d\to Y\) is finite. This gives \(H^2(K_d,\QQ)\ne0\). The Hopf exact sequence gives a surjection \(H_2(X_d,\ZZ)\to H_2(K_d,\ZZ)\), and \(H^2(K_d,\QQ)=\operatorname{Hom}(H_2(K_d,\ZZ),\QQ)\) by the universal coefficient theorem. Thus \(H_2(K_d,\ZZ)\) is finitely generated of positive rank. The same theorem gives \(H^2(K_d,\mathbb F_p)\ne0\) for every prime \(p\). The continuous degree-two cohomology of \(\widehat K\) is zero because \(\widehat K=1\). This proves Corollary~\ref{cor:applications}.

\begin{corollary}\label{cor:wedge}
The universal cover of each surface in Theorem~\ref{thm:main} has the homotopy type of a countably infinite wedge of two-spheres.
\end{corollary}
\begin{proof}
By Stein Morse theory \cite{AF,CE}, the Stein universal cover \(\widetilde X\) has the homotopy type of a countable CW complex of dimension at most two. Indeed, the real Hessian of a strictly plurisubharmonic Morse exhaustion cannot be negative on both \(v\) and \(Jv\), so every Morse index is at most the complex dimension. A basis of the free abelian group \(H_2(\widetilde X,\ZZ)\) is represented by spheres by the Hurewicz theorem. The resulting map from their wedge is a homotopy equivalence by the homological Whitehead theorem.

If this group had finite rank, the action of \(K\) on \(H^2(\widetilde X,\ZZ)\) would be trivial by applying Corollary~\ref{cor:applications} after tensoring with \(\QQ\). Choose a subgroup \(C_p\subset K\) of prime order, which exists by Theorem~\ref{thm:main}. The integral cohomology spectral sequence of the Borel fibration over \(BC_p\) abuts to \(H^*(\widetilde X/C_p,\ZZ)\), since the action is free \cite[Th.~5.2, p.~135]{McCleary}. Its two possible rows \(H^i(C_p,H^j(\widetilde X,\ZZ))\), \(j=0,2\), are supported in even columns by the cyclic group resolution \cite[Ch.~I, (6.3)]{Brown}. The only possible differential \(d_3\) changes column parity. Hence the bottom-row groups \(H^{2m}(C_p,\ZZ)=\ZZ/p\) for \(m>0\) survive, contradicting the real dimension four of the quotient. The rank is therefore countably infinite.
\end{proof}

Dimca, Papadima and Suciu constructed projective examples with Stein universal cover of infinite wedge type \cite[Th.~A and Cor.~5.4]{DPS}. Here the wedge type occurs together with trivial profinite completion and a prescribed finite subgroup spectrum.

In each even ambient dimension \(n\ge4\), the root stacks in \cite[Prop.~1.10]{Eys} are smooth, effective and proper with projective coarse spaces. Applying Lemma~\ref{lem:realization} in dimension \(n\) realizes their groups, for all root orders, as fundamental groups of smooth projective surfaces with ample canonical bundle. This answers \cite[Question~1.11]{Eys}.

\end{document}